\documentclass[10pt,a4paper,leqno]{article}

\usepackage[T1]{fontenc}
\usepackage[latin1]{inputenc}
\usepackage[francais]{babel}
\usepackage{amsmath}
\usepackage{amssymb}
\usepackage{graphicx}

\newtheorem{definition}{Définition}

\newtheorem{theorem}{Théorème}     
\newtheorem{lem}{Lemme}
\newtheorem{coro}{Corollaire}

\newtheorem{rem}{Remarque}
\newenvironment{dem}{\smallskip\noindent{\bf Démonstration}~:%
{\nopagebreak[0]}}%
{\nopagebreak[0]\hfill$\Box$ \medskip}%

\newcommand{\R}{\mbox{$\mathbb{R}$}}

\newcommand{\ben}{{\rm ben\,}}
\newcommand{\benp}{{\rm ben}_p\,}

\def\ds{\displaystyle}
\def\iy{\infty}

\def\al{\alpha}
\def\be{\beta}
\def\ga{\gamma}

\def\de{\delta}

\def\la{\lambda}

\def\te{\theta}

\def\vep{\varepsilon}
\def\vfi{\varphi}
\def\om{\omega}

\def\si{\sigma}

\def\ce{\mathcal{E}}

\def\cn{\mathcal{N}}

\def\ni{\noindent}
\def\mk{\medskip}
\def\bk{\bigskip}
\def\ap{\approx}

\def \dv{\,|\,}

\def\Lr{\Longrightarrow}

\title{Quelques inégalités effectives entre des fonctions arithmétiques 
usuelles}
\author{Jean-Louis NICOLAS\,\footnote{Recherche financ\'ee
par le CNRS, Institut Camille Jordan, UMR 5208}.}

\begin{document}
\maketitle

\hfill
\begin{minipage}[t]{72mm}
A Wladyslaw Narkiewicz pour son soixante-dixième 
anniversaire, en très amical hommage.
\end{minipage}
\bigskip

\medskip
\ni
{\bf Abstract.} Let us denote by $\tau(n)$ and
$\si(n)$ the number and the sum of the divisors of $n$ and by $\vfi$ 
Euler's function. We give effective 
upper bounds for $\frac{n}{\vfi(n)}$ in terms of $\vfi(n)$, and for 
$\frac{\si(n)}{n}$ in terms of $\tau(n)$. 

\section{Introduction}

Soit $n$ un entier positif. Nous utilisons les fonctions arithmétiques 
classiques:
$$\tau(n)=\sum_{d\dv n} 1, \quad\si(n)=\sum_{d\dv n} 1, \quad
\om(n)=\sum_{\substack{p\dv n\\ p\;\text{ premier }}} 1, \quad
\pi(x)=\sum_{\substack{p\;\le\; x\\ p\;\text{ premier }}} 1$$
tandis que $\vfi$ désigne la fonction d'Euler.
On note $p_k$ le $k$-ième nombre premier et $\ga\ap 0.57721566$ 
la constante d'Euler.

Dans cet article, nous nous intéressons aux grandes valeurs des   
fonctions $\frac{\si(n)}{n}$  et $\frac{n}{\vfi(n)}$. L'ordre maximum de ces 
deux fonctions est $e^\ga \log \log n$ (cf. \cite[Th. 323 et Th. 328]{HW}).

De façon plus précise, Rosser et Schoenfeld ont montré dans \cite{RS} 
\begin{equation}\label{RSphi}
\frac{\si(n)}{n}\le \frac{n}{\vfi(n)}\le e^\ga \log \log n+
\frac{2.50637}{\log \log n},\quad n\ge 3
\end{equation} 
tandis qu'il est prouvé dans \cite{Njtn} qu'il existe une infinité 
de nombres $n$ pour lesquels $\frac{n}{\vfi(n)} > e^\ga \log \log n$.
Notons que $e^\ga\ap 1.7810724$

Le comportement de $\frac{\si(n)}{n}$ est différent: dans \cite{Rob}, 
il est démontré que l'hypothèse de Riemann est équivalente à
$$ \forall n \ge 5041, \qquad \frac{\si(n)}{n}\le  e^\ga \log \log n.$$

Nous prouverons

\begin{theorem}\label{th1}
Pour $n\ge 3$, on a
\begin{equation}\label{3.65}
\frac{n}{\vfi(n)}\le e^\ga \log \log \vfi(n)+3.65278\ldots
\end{equation}
avec égalité pour $n=6$; pour $n\ge 211$, on a
\begin{equation}\label{3}
\frac{n}{\vfi(n)}\le 3\log \log \vfi(n)
\end{equation}
et pour $n\ge 7$, on a 
\begin{equation}\label{2.95}
\frac{n}{\vfi(n)}\le e^\ga \log \log \vfi(n)+ 
\frac{2.95503\ldots}{\log \log \vfi(n)}
\end{equation}
avec égalité pour $n=30030=2\times 3\times 5\times 7\times 11\times 13$.
\end{theorem}
L'in\'egalit\'e \eqref{3.65} r\'epond \`a une question pos\'ee par
A. Schinzel dans \cite{Sch}.

\mk

Dans l'article \cite[formule (17)]{Ten1}, G. Tenenbaum d\'emontre la relation
\begin{equation}\label{tenb}
\frac{\si(n)}{n}\ll \log\log (2\tau(n)).
\end{equation}
Le th\'eor\`eme suivant donne une forme effective \`a cette in\'egalit\'e.
\begin{theorem}\label{th2}
On a pour tout $n\ge 2$
\begin{equation}\label{in3}
\frac{\si(n)}{n}\le 2.59790\ldots \log\log \;(3\tau(n))
\end{equation}
avec égalité pour $\ds n=M_1\stackrel{def}{=\!=}
 2^83^55^3 7^211^213^2 \prod_{17\le p\le 113} p$, et
\begin{equation}\label{in4}
\frac{\si(n)}{n}\le e^\ga(\log\log\; (e\tau(n))+\log\log\log\; (e^e \tau(n))
+0.941444079\ldots
\end{equation}
avec égalité pour
$$n=M_2 \stackrel{def}{=\!=} 2^{18}3^{11}5^77^6
\prod_{11\le p \le 19}\; p^4\prod_{23\le p \le 47}\; p^3
\prod_{53\le p \le 277}\; p^2 \prod_{281\le p \le 45439}\; p.$$ 
\end{theorem}

Les coefficients $3$ dans \eqref{in3} et $e$ et $e^e$ dans \eqref{in4} 
pourraient être modifiés. Ces quantités ont été choisies assez grandes de
façon que les valeurs prises par la fonction $\log\log$ pour les petites 
valeurs de la variable ne soient pas prépondérantes. 

La démonstration des formules des théorèmes \ref{th1} et \ref{th2} 
se fait de la façon 
suivante: on commence par les prouver pour $n$ ou $\tau(n)$
suffisamment grand. Il reste alors un nombre fini de valeurs de 
$n$ ou $\tau(n)$ à examiner. Mais ce nombre de valeurs est très grand, 
ce qui interdit une étude systématique par 
ordinateur. On construit alors une sous-famille beaucoup plus petite de 
nombres pour lesquels il suffira de faire les calculs.

Pour le théorème \ref{th1}, cette sous-famille est constituée par les nombres 
$N_k=\prod_{1\le i\le k} p_i$.  Pour le théorème \ref{th2}, ce sont 
les nombres $(\si,\tau)$--superchampions définis au paragraphe \ref{super} et 
qui ressemblent aux nombres hautement composés supérieurs  introduits par 
Ramanujan (cf. \cite{Ram,Rampos}).

Dans l'article \cite{RS}, Rosser et Schoenfeld démontrent un peu plus que
la formule \eqref{RSphi}. En fait, ils montrent que cette formule 
dans laquelle la constante $2.50637$ est remplacée par $5/2$
est vérifiée pour tout nombre $n\ge 3$ à l'exception de
$n=223092870=2\times3\times5\times7\times11\times13\times17\times19
\times23$. On peut obtenir un résultat similaire pour les majorations des 
théorèmes \ref{th1} et \ref{th2}. Nous le ferons explicitement pour 
\eqref{in3}.

\begin{theorem}\label{th3}
L'inégalité
\begin{equation}\label{in3bis}
\frac{\si(n)}{n}\le \frac{2597}{1000}\log\log \;(3\tau(n))
\end{equation}
est vérifiée pour tout $n\ge 2$ à l'exception de 12 nombres.
\end{theorem}

La démonstration du théorème \ref{th3}  sera donnée 
au paragraphe \ref{par5}.
Elle utilise la notion de bénéfice qui précise le comportement d'un nombre 
ordinaire par rapport à un nombre $(\si,\tau)$--superchampion.

\section{Démonstration du théorème \ref{th1}}

Nous notons $N_k=2\times 3\times \ldots \times p_k$ 
le produit des $k$ premiers nombres premiers. 
Démontrons d'abord trois lemmes.

\subsection{Trois lemmes}

\begin{lem}\label{lem1}
Soit $k\ge 1$ et un entier $n$ vérifiant $N_k < n < N_{k+1}$. Alors on a 
\begin{equation}\label{nNk}
\frac{n}{\vfi(n)} < \frac{N_k}{\vfi(N_k)}=
\prod_{2\;\le \;p\; \le \;p_k} \frac{p}{p-1}\cdot
\end{equation}
Soit $\ell\ge 1$ et $n > N_\ell$; on a
\begin{equation}\label{fiNl}
\vfi(n) > \vfi(N_\ell)=\prod_{2\;\le \;p\; \le \;p_\ell} (p-1).
\end{equation}
\end{lem}

\begin{dem}
Soit $n=q_1^{\al_1}q_2^{\al_2}\ldots q_j^{\al_j}$ la décomposition 
en facteurs premiers de $n$.
Puisque $n < N_{k+1}$, le nombre $j=\om(n)$ vérifie $j\le k$ et l'on a
\begin{equation}\label{maj1}
\frac{n}{\vfi(n)}=\prod_{i=1}^j \frac{1}{1-1/q_i}\le 
\prod_{i=1}^j \frac{1}{1-1/p_i}=\frac{N_j}{\vfi(N_j)} 
\le \frac{N_k}{\vfi(N_k)}\cdot 
\end{equation}  
De plus, la première inégalité de \eqref{maj1} est stricte 
car $n\ne N_k$ et l'on 
ne peut avoir $j=k$ et $q_i=p_i$ pour $1\le i \le k$. Ceci démontre la 
relation \eqref{nNk}. 

Prouvons maintenant \eqref{fiNl}. Puisque $n > N_\ell$, il existe $k\ge \ell$
tel que $N_k\le n < N_{k+1}$. 

Si $k=\ell$, on a $N_\ell=N_k < n < N_{k+1}$ et \eqref{nNk} entra\^{i}ne 
$\frac{n}{\vfi(n)} < \frac{N_\ell}{\vfi(N_\ell)}$,
soit $\vfi(n) > \frac{n}{N_\ell}\vfi(N_\ell) > \vfi(N_\ell)$.

Si $k > \ell$, \eqref{nNk} entraîne $\frac{n}{\vfi(n)} \le 
\frac{N_k}{\vfi(N_k)}$, et l'on a 
$\vfi(n) \ge \frac{n}{N_k}\vfi(N_k) \ge \vfi(N_k) > \vfi(N_\ell)$.
\end{dem}

\begin{rem}
Les nombres $M$ tels que $m > M\;\;\Lr \;\; \vfi(m) > \vfi(M)$ ont été appelés
``sparsely totient'' et étudiés par Masser et Shiu (cf. \cite{MS}). La formule
\eqref{fiNl} montre que les nombres $N_\ell$ sont ``sparsely totient''.
\end{rem}


\begin{lem}\label{lem2}
Soit $k\ge 1$ et $f: [\vfi(N_k),+\iy)\to \R$ une fonction strictement 
croissante vérifiant 
$$\frac{N_k}{\vfi(N_k)}\le f(\vfi(N_k)).$$
Soit $n$ vérifiant $N_k < n < N_{k+1}$, on a 
$$\frac{n}{\vfi(n)} < f(\vfi(n)).$$ 
\end{lem}

\begin{dem}
Par le lemme \ref{lem1} et nos hypothèses, on a $\vfi(n) > \vfi(N_k)$ et 
$$\frac{n}{\vfi(n)} < \frac{N_k}{\vfi(N_k)}\le f(\vfi(N_k)) < f(\vfi(n)).$$
\end{dem}

\begin{lem}\label{lem3}
Pour $n\ge 7$, on a $\vfi(n)\ge \sqrt n$.
\end{lem}

\begin{dem}
Utilisons la minoration très simple (cf. \cite[p. 319 ]{Rib}) :
$$\vfi(n)\ge \frac{n\log 2}{\log (2n)}\cdot$$
Il est facile de voir que pour $t$ réel, $t\ge 40$, on a 
$\frac{t\log 2}{\log (2t)}\ge \sqrt t$, puis on vérifie le lemme 
pour $7\le n\le 39$.
\end{dem}

\subsection{Les grandes valeurs de $\mathbf{n}$ : $\mathbf{n\ge N_{14}}$}

Posons $b:=2.51$ et $g(t)=e^\ga\log\log t +\frac{b}{\log\log t}$. La 
fonction $g$ est croissante pour $t \ge 27 > \exp(\exp(\sqrt{be^{-\ga}}))$.

Soit $n\ge N_6=30030$. Le lemme \ref{lem1} entraîne $\vfi(n)\ge \vfi(N_6)=5760$
tandis que, par le lemme \ref{lem3}, on a
$$\log\log n\le \log\log \vfi(n)+\log 2.$$
Par \eqref{RSphi} et la croissance de $g$, il suit
\begin{eqnarray}\label{2.12}
\frac{n}{\vfi(n)}\le g(n) &\le& e^\ga \left(\log\log \vfi(n)+\log 2\right)
+\frac{b}{\log\log \vfi(n)+\log 2}\notag\\
&\le& e^\ga\log\log \vfi(n)+e^\ga\log 2+\frac{b}{\log\log 5760+\log 2}\notag\\
&\le& e^\ga\log\log \vfi(n)+ 2.12\;.
\end{eqnarray}
On déduit de l'inégalité \eqref{2.12}
$$\frac{n}{\vfi(n)}\le 
\log\log \vfi(n) \left[e^\ga+\frac{2.12}{\log\log \vfi(n)}\right]
\le \log\log \vfi(n) \left[e^\ga+\frac{2.12}{\log\log 5760}\right]$$
ce qui montre, pour $n\ge 30030$,
\begin{equation}\label{3fi1}
\frac{n}{\vfi(n)}\le 3\log\log \vfi(n).
\end{equation}
Pour $n\ge 30030$, il vient alors en posant
\begin{eqnarray*}
u &=&\log\log\vfi(n)\ge \log\log 5760 > 2.15,\\
\log n &\le& \log \vfi(n) \left(1+ \frac{\log\log\log\vfi(n)+\log 3}
{\log \vfi(n)}\right)\\
\log \log n &\le& \log\log \vfi(n) + \frac{\log\log\log\vfi(n)+\log 3}
{\log \vfi(n)}=u+\frac{h(u)}{u}
\end{eqnarray*}
avec
$$h(u)=u(\log u+\log 3)e^{-u}.$$
Par la croissance de la fonction $g$ et \eqref{RSphi}, 
la majoration ci-dessus implique
\begin{equation}\label{nhu}
\frac{n}{\vfi(n)}\le 
e^\ga\left(u+\frac{h(u)}{u}\right)+\frac{b}{u+h(u)/u}\le
e^\ga u + \frac{b+e^\ga h(u)}{u}\cdot
\end{equation} 

Or la fonction $u\mapsto h(u)$ est décroissante pour $u> 1.64$, et l'on a 
$e^\ga h(u) < 0.43$ pour $u \ge 3.55$.

Pour $n\ge N_{14}=\prod_{p\le 43} \; p$, on a par le lemme \ref{lem1}, 
$\vfi(n)\ge \vfi(N_{14})=1.85\ldots\;10^{15}$, $u=\log\log \vfi(n)\ge 
\log\log \vfi(N_{14}) > 3.55$ et \eqref{nhu} implique
\begin{equation}\label{nN14}
\frac{n}{\vfi(n)}\le e^\ga \log \log \vfi(n)+ 
\frac{2.94}{\log \log \vfi(n)}, \quad n\ge N_{14}.
\end{equation}

\subsection{Les petites valeurs de $n$}

\subsubsection{Démonstration de \eqref{3.65}}

La formule \eqref{2.12} a été établie pour $n\ge 30030$. Pour 
démontrer \eqref{3.65}, il suffit de calculer 
$G(n)=\frac{n}{\vfi(n)}-e^\ga \log\log\vfi(n)$ pour $3\le n\le 30029$.

Mais on peut se contenter de calculer
$$
\begin{array}{|r|ccccccc|}
n=&3&4&5&N_2=6&N_3=30&N_4=210&N_5=2310\\
\hline
G(n)=&2.15&2.65&0.67&3.65&2.45&1.96&1.57\\
\hline
\end{array}
$$
et d'appliquer le lemme \ref{lem2} à la fonction
$$f(t)=e^\ga\left(\log\log t-\log\log \vfi(6)\right)+\frac{6}{\vfi(6)}
=e^\ga \log\log t+G(6).$$
Pour $k\in\{2,3,4,5\}$, on a $\frac{N_k}{\vfi(N_k)}-f(\vfi(N_k))=
G(N_k)-G(6)\le 0$ et il en résulte 
$\frac{n}{\vfi(n)}\le f(\vfi(n))=e^\ga\log\log \vfi(n)+G(6)$ 
pour $N_2=6\le n < N_6=30030$.

\subsubsection{Démonstration de \eqref{3}}

La formule \eqref{3fi1} prouve \eqref{3} pour $n\ge N_6$.
Comme $\frac{N_5}{\vfi(N_5)\log\log\vfi(N_5)}=2.64\ldots$, le 
lemme \ref{lem2} appliqué à la fonction $f(t)=3\log\log t$ démontre 
\eqref{3} pour $N_5=2310\le n < N_6=30030$.

Il reste à calculer $\frac{n}{\vfi(n)}-3\;\log\log\vfi(n)$ pour $n\le 2309$.
Cette quantité est non définie pour $n\in \{1,2\}$ et positive pour
$$n\in\{3,4,5,6,8,10,12,14,18,20,24,30,36,42,60,66,84,90,120,210\}.$$
Pour $n\ge 3$, le maximum de $\frac{n}{\vfi(n)\log \log \vfi(n)}$ est 
atteint pour $n=12$ et vaut $9.18458\ldots$

\subsubsection{Démonstration de \eqref{2.95}}

Posons $\ds c_k=\log\log\vfi(N_k)\left[\frac{N_k}{\vfi(N_k)}-e^\ga
\log\log\vfi(N_k)\right]$.
On calcule
$$
\begin{array}{|r|ccccccccccc|}
k=&4&5&6&7&8&9&10&11&12&13&14\\
\hline
c_k=&2.66&2.86&2.96&2.92&2.94&2.93&2.82&2.77&2.68&2.59&2.55\\
\hline
\end{array}
$$
et l'on pose

$$c=\max_{4\,\le\,k\,\le\,14} c_k=c_6\ap 2.9550377$$
La fonction $f(t)=e^\ga\log\log t+\frac{c}{\log\log t}$ est croissante pour 
$t\ge 38 > e^{e^{\sqrt{ce^{-\ga}}}}$. Comme $\vfi(N_4)=\vfi(210)=48$, on 
applique le lemme \ref{lem2} à la fonction $f$ pour $4\le k\le 14$; 
cela prouve \eqref{2.95} pour $N_4\le n < N_{15}$ et, compte tenu de 
\eqref{nN14}, pour $n\ge N_{4}=210$.

Il reste à vérifier \eqref{2.95} pour $n\le 209$. Les exceptions sont
$n\in\{1,2,3,4,6\}$; Pour ces valeurs de $n$, on a $\vfi(n)=1$ ou $2$ et
$\log\log \vfi(n)$ est soit non défini soit négatif. 

\section{Les nombres $\mathbf{(\si,\tau)}$--superchampions}\label{super}

Soit $\vep$ un nombre réél positif. Il résulte de \eqref{tenb} que
$\lim_{n\to\iy} \frac{\si(n)}{n\tau(n)^\vep}=0$ et que 
$\frac{\si(n)}{n\tau(n)^\vep}$ est borné.

\begin{definition}\label{defsupch}
On dit que $N$ est un nombre $(\si,\tau)$--superchampion s'il existe $\vep > 0$
tel que, pour tout $n\ge 1$, on ait 
\begin{equation}\label{ineqsupch}
\frac{\si(n)}{n\tau(n)^\vep} \le \frac{\si(N)}{N\tau(N)^\vep}\cdot
\end{equation}
\end{definition}
Cette définition suit l'exemple de Ramanujan qui dans \cite[\S 32]{Ram} a 
introduit les nombres hautement composés supérieurs. Notons qu'elle 
entraîne la propriété suivante: 
\begin{equation}\label{champ}
\tau(n) \le \tau(N)\quad \Lr \quad \frac{\si(n)}{n}\le \frac{\si(N)}{N}.
\end{equation}
En effet, par la définition \ref{defsupch}, on a pour tout $n\ge 1$ vérifiant
$\tau(n)\le \tau(N)$
$$\frac{\si(n)}{n} \le \frac{\si(N)}{N}
\left(\frac{\tau(n)}{\tau(N)}\right)^\vep\le \frac{\si(N)}{N}\cdot$$

La fonction $\frac{\si(n)}{n\tau(n)^\vep}$ est multiplicative; pour 
trouver son maximum, on recherche d'abord pour chaque $p$ premier le 
maximum de la fonction 
$\al\mapsto \frac{\si(p^\al)}{p^\al \tau(p^\al)^\vep}\cdot$

\begin{lem}\label{lempsi}
Posons pour $p$ premier et $\al$ entier, $\al\ge 1$,
\begin{equation}\label{psi}
\psi(p,\al)=
\frac{\log \left(1+\frac{1}{p+p^2+\ldots+p^\al}\right)}
{\log\left(1+\frac1\al\right)}=
\frac{\log \left(1+\frac{p-1}{p^{\al+1}-p}\right)}
{\log\left(1+\frac1\al\right)}=
\frac{\log \frac{p^{\al+1}-1}{p^{\al+1}-p}}
{\log\left(1+\frac1\al\right)}
\end{equation}
et $\psi(p,0)=+\iy$. La fonction $\psi(p,\al)$ est strictement
décroissante en $p$ et en $\al$. On a 
$$\lim_{p\to\iy} \psi(p,\al)=0\qquad \text{ et } \qquad 
\lim_{\al\to\iy} \psi(p,\al)=0.$$
Soit $\al_p$ un nombre en lequel la fonction 
$\al\mapsto \frac{\si(p^\al)}{p^\al \tau(p^\al)^\vep}$ atteint son maximum.
Alors, on a 
\begin{equation}\label{psieps}
\psi(p,\al_p+1) \le \vep \le \psi(p,\al_p).
\end{equation} 
\end{lem}

\begin{dem}
La décroissance en $p$ et les limites sont faciles à établir. Pour 
prouver la décroissance en $\al$, nous utiliserons les inégalités
$$\frac{1}{1+t} < \log\left(1+\frac 1t\right) < \frac 1t, \qquad t > 0.$$
On a $\psi(p,1) < \psi(p,0)=+\iy$. Supposons maintenant $\al\ge 2$; il vient
\begin{eqnarray*}
\psi(p,\al) &<& \frac{1}{(p+p^2+\ldots+p^\al)\log(1+1/\al)}\\
&=& \frac{1}{p \log(1+1/\al)}\;
\frac{1}{1+p+\ldots+p^{\al-1}}\\
&<& \frac{1}{p\log(1+1/\al)}
\log\left(1+\frac{1}{p+p^2+\ldots+p^{\al-1}}\right)\\
&=& \psi(p,\al-1)\frac{\log\left(1+\frac{1}{\al-1}\right)}
{p\log\left(1+\frac 1\al\right)}\le
\psi(p,\al-1)\frac{\log\left(1+\frac{1}{\al-1}\right)}
{2\log\left(1+\frac 1\al\right)}
\end{eqnarray*}
et l'on a 
$$\left(1+\frac 1\al\right)^2> 1+\frac 2\al \ge 1+\frac{1}{\al-1}$$
ce qui achève la preuve de la décroissance de $\psi(p,\al)$ en $\al$.

On pose ensuite $\te(p,0)=0$ et, pour $\al\ge 1$,

\begin{equation}\label{te}
\te(p,\al)=\log\frac{\si(p^\al)}{p^\al}-\vep\log \tau(p^\al)
=\log\left(1+\frac{1}{p}+\ldots+\frac{1}{p^\al}\right)
-\vep\log(\al+1).
\end{equation}
Pour $\al\ge 1$, on a
\begin{equation}\label{tem}
\te(p,\al-1)-\te(p,\al)=\log\left(1+\frac 1\al\right)(\vep-\psi(p,\al))
\end{equation}
et \eqref{psieps} en résulte compte tenu de la décroissance 
en $\al$ de la fonction $\psi$.
\end{dem}
\begin{lem}\label{lemexpo}
Soit $N$ un nombre $(\si,\tau)$--superchampion de paramètre $\vep$ 
(autrement dit la fonction $n\mapsto \frac{\si(n)}{n\tau(n)^\vep}$ 
atteint son maximum en $N$). Un diviseur premier $p$ de $N$ vérifie
$p\le \frac{1}{2^\vep-1}$. Si $p > q$ sont deux nombres premiers, on a 
$$v_p(N) \le v_q(N).$$
\end{lem}

\begin{dem}
Par le lemme \ref{lempsi}, le nombre $N$ s'écrit $N=\prod_{p} \;p^{\al_p}$ 
où, pour chaque nombre premier $p$, $\al_p$ vérifie \eqref{psieps}.

L'inégalité $p > \frac{1}{2^\vep-1}$ est équivalente à 
$\vep > \frac{\log(1+1/p)}{\log 2}=\psi(p,1)$, ce qui, par \eqref{psieps},
implique $\al_p=0$.

Par la décroissance en $p$ de $\psi(p,\al)$ annoncée dans le 
lemme \ref{lempsi} et par \eqref{psieps}, on a
$$\psi(p,\al_p+1)\le \vep\le \psi(p,\al_p)< \psi(q,\al_p)$$
ce qui, encore par \eqref{psieps}, entraîne $\al_q\ge \al_p$.
\end{dem}

\subsection{Détermination des nombres $\mathbf{(\si,\tau)}$--superchampions}

Posons
\begin{equation}\label{E}
\ce_p=\{\psi(p,\al), \al=1,2,\ldots\}\quad \text{ et } \quad 
\ce=\left(\bigcup_{p\; \text{ premier}} \ce_p\right) \cup \{+\iy\}.
\end{equation}
D'après le théorème des six exponentielles (cf. \cite[p. 14]{Wal}), si 
$p,q,r$ sont trois nombres premiers distincts, les ensembles $\ce_p,\ce_q$
et $\ce_r$ ont une intersection vide (cf. \cite[p. 455]{AlaErd}
et \cite[p. 71]{ErdNic}).

En fait, si $p$ et $q$ sont deux nombres premiers distincts,  
il est vraisemblable que l'on a $\ce_p\cap \ce_q=\emptyset$.  
Dans les calculs effectués dans cet article, tous les nombres 
$\vep \ge \psi(2248723,1)\ap0.00000064156
$
et vérifiant $\vep\in\ce$ appartiennent à un seul ensemble $\ce_p$.

Pour mesurer la proximité de deux éléments distincts $\vep'$ et $\vep''$ de
$\ce$, il est commode de calculer $|g_1(\vep')-g_1(\vep'')|$ avec 
$g_1(\vep)=\frac{1}{2^\vep-1}$. Notons que l'on a $g_1(\psi(p,1))=p$.
Pour $\vep'' > \vep' \ge \psi(2248723,1)$, nous avons trouvé
$$\min(g_1(\vep')-g_1(\vep''))=
g_1(\psi(71453,1))-g_1(\psi(349,2))\ap 0.0381\;.$$

Ordonnons les éléments de l'ensemble $\ce$:
$$\ce=\{\vep_0=+\iy > \vep_1 > \vep_2 > \ldots \}.$$
Le lemme suivant, qui est voisin de la proposition 4 de \cite{ErdNic},
détermine les nombres $(\si,\tau)$--superchampions:
\begin{lem}\label{lemsupch}
(i) Soit $i\ge 0$ et $\vep_{i} > \vep > \vep_{i+1}$. Le maximum de la fonction 
$n\mapsto\frac{\si(n)}{n\tau(n)^\vep}$ est atteint en un seul nombre
$$N^{(i)}=\prod_{p < \frac{1}{2^\vep-1}} \; p^{\al_p}$$
avec $\al_p$ défini de façon unique par \eqref{psieps}.

(ii) Si $\vep=\vep_i$ et si $\vep$ appartient à un seul ensemble $\ce_p$, 
la fonction $n\mapsto\frac{\si(n)}{n\tau(n)^\vep}$ atteint 
son maximum en deux points $N^{(i-1)}$ et $N^{(i)}=pN^{(i-1)}$.

(iii) Si $\vep=\vep_i$ et si $\vep\in \ce_q\cap\ce_r$, 
la fonction $n\mapsto\frac{\si(n)}{n\tau(n)^\vep}$ atteint 

son maximum en quatre points $N^{(i-1)}, qN^{(i-1)},rN^{(i-1)}$ et 
$N^{(i)}=qrN^{(i-1)}$.
\end{lem}

\begin{dem}
La démonstration s'appuie sur les lemmes \ref{lempsi} et \ref{lemexpo}. 
On observera que si $\vep=\vep_i=\psi(p,\al)$, il y a 
dans \eqref{psieps} deux valeurs possibles pour $\al_p$, $\al_p=\al$ 
et $\al_p=\al-1$.
\end{dem}

$$
\begin{array}{|l|r|l|r|c|}
i&N^{(i)}&\frac{\si(N^{(i)})}{N^{(i)}}&\tau(N^{(i)})&\vep\\
\hline
0&1&1&1&\vep_0=+\iy > \vep \ge \vep_1=\log(3/2)/\log 2=0.585\\
1&2&1.5&2&\vep_1\ge \vep \ge \vep_2=\log(4/3)/\log 2=0.415\\
2&6&2&4&\vep_2\ge \vep\ge \vep_3=\log(7/6)/\log(3/2)=0.380\\
3&12&2.333&6&\vep_3\ge \vep\ge \vep_4=\log(6/5)/\log 2=0.263\\
4&60&2.8&12&\vep_4\ge \vep\ge \vep_5=\log(15/14)/\log(4/3)=0.240\\
5&120&3&16&\vep_5\ge \vep\ge \vep_6=\log(13/12)/\log(3/2)=0.197\\
6&360&3.25&24&\vep_6\ge \vep\ge \vep7=\log(8/7)/\log 2=0.193\\
7&2520&3.714&48&\vep_7\ge \vep\ge \vep_8=\log(31/30)/\log(5/4)=0.147\\
8&5040&3.838&60&\vep_8\ge \vep\ge \vep_9=\log(12/11)/\log 2=0.126\\
9&55440&4.187&120&\vep_9\ge \vep\ge \vep_{10}=\log(14/13)/\log 2=0.107\\
\hline
\end{array}
$$
\medskip

On lit dans la table ci-dessus que les deux nombres $2$ et $120$ sont 
$(\si,\tau)$--super\-champions pour les paramètres $1/2$ et $1/5$ 
respectivement; de la relation \eqref{ineqsupch}, on déduit alors  
les inégalités valables pour $n\ge 1$, 
\begin{equation}\label{in1}
\frac{\si(n)}{n}\le \frac{\si(2)}{2}\sqrt{\frac{\tau(n)}{\tau(2)}}=
\frac{3}{2\sqrt 2} \sqrt{\tau(n)} \le 1.061\;\sqrt{\tau(n)}
\end{equation}
et
\begin{equation}\label{in2}
\frac{\si(n)}{n}\le 
\frac{\si(120)}{120}\left(\frac{\tau(n)}{\tau(120)}\right)^{0.2}= 
3\left(\frac{\tau(n)}{16}\right)^{0.2}\le 
1.72305\; \tau(n)^{0.2}.
\end{equation}

\begin{definition}\label{defN+}
Soit $\vep > 0$. On note $N^+_\vep$ (resp. $N^-_\vep$) 
le plus grand (resp. petit) nombre $(\si,\tau)$--superchampion de paramètre
$\vep$. 
\end{definition} 

\begin{rem}\label{remN+}
D'après le lemme \eqref{lemsupch}, si $\vep\notin\ce$
(cas (i)), on a $N^+_\vep =N^-_\vep$. Si $\vep=\vep_i\in\ce$ 
(cas (ii) et cas (iii)), on a $N^+_\vep=N^{(i)}$ et $N^-_\vep=N^{(i-1)}$.

Le nombre $N^{(i)}$ est un diviseur de $N^{(i+1)}$. Si un nombre
$(\si,\tau)$--supercham\-pion $N$ vérifie $N\ge N^+_\vep$, alors, $N$ est un 
multiple de $N^+_\vep$.
\end{rem}

\subsection{Aspect géométrique}

A chaque entier $n\ge 1$, associons dans un système d'axes $xOy$ le point 
$(\log \tau(n),\log (\si(n)/n))$ appelé {\it image} de $n$. La droite 
$D(\vep,n)$ de pente $\vep$
et passant par l'image de $n$ coupe l'axe $Oy$ au point d'ordonnée
$\log \frac{\si(n)}{n}-\vep\log\tau(n)$.

Il résulte des calculs précédents que l'ensemble des images des entiers 
$n\ge 1$ a une enveloppe convexe qui est une ligne polygonale dont les pentes 
des côtés sont les éléments de $\ce$ et dont les sommets sont les 
images des nombres $(\si,\tau)$--superchampions (cf. \textsc{Fig}. 1).

Soit $\vep > 0$ fixé. 

Si $\vep \notin\ce$, la plus haute droite 
$D(\vep,n)$ va passer par un seul des sommets de l'enveloppe convexe.

Si $\vep\in\ce$ et si $\vep$ n'appartient qu'à un seul $\ce_p$ (cas (ii) du 
lemme \ref{lemsupch}), la plus haute droite 
$D(\vep,n)$ va joindre deux sommets consécutifs de l'enveloppe convexe. 

Si $\vep\in\ce_q\cap\ce_r$ (cas (iii) très peu probable du 
lemme \ref{lemsupch}), la plus haute droite 
$D(\vep,n)$ contiendra les images de quatre nombres 
$(\si,\tau)$--superchampions.

\hspace{-6cm}\includegraphics[scale=1.]{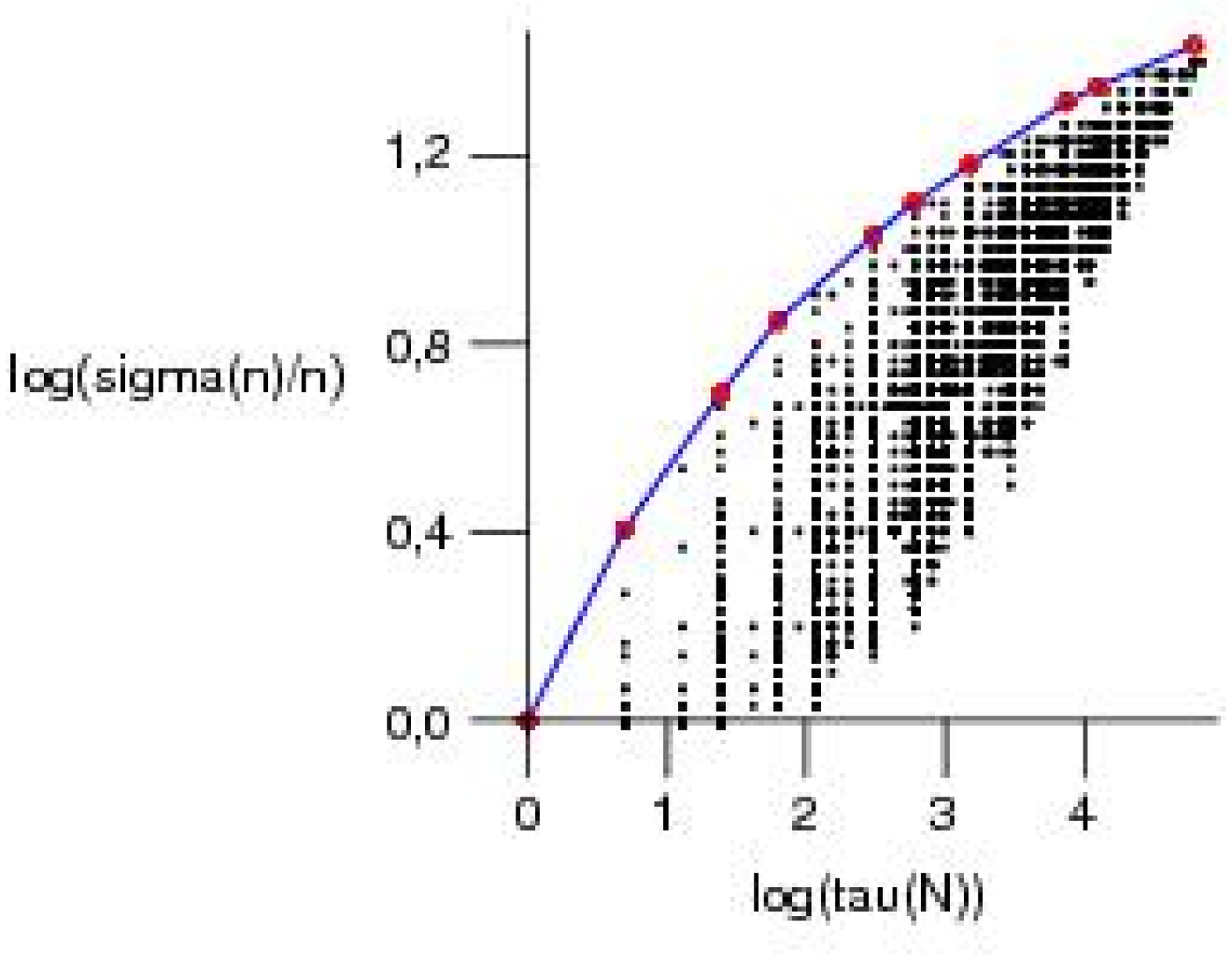}  

\vspace{-25mm}
 
\begin{center}

\textsc{Fig}. 1 : Les images des nombres $1\le n\le 55440$.

\end{center}

\section{Démonstration du théorème \ref{th2}}

\subsection{Deux lemmes}
 
\begin{lem}\label{lem4}
Soit $F\!:\! [u_0,u_1]\subset \R \to \R_+^*$ une fonction logarithmiquement 
concave. Soit $\vep\in\ce$ (cf. \eqref{E}) et 
$N'$ et $N''$ deux nombres $(\si,\tau)$--superchampions
de même paramètre $\vep$ (cf. lemme \ref{lemsupch}). 
Soit $n$ un nombre entier vérifiant
$$e^{u_0}\le \tau(N')\le \tau(n)\le \tau(N'')\le e^{u_1}.$$
Alors, on a
$$\frac{\si(n)}{n\;F(\log(\tau(n)))}\le \max
\left(\frac{\si(N')}{N'\;F(\log(\tau(N')))}\;,\;\;
\frac{\si(N'')}{N''\;F(\log(\tau(N'')))}\right).$$
\end{lem}

\begin{dem}
La démonstration s'inspire de celle de la proposition 1 de \cite{Rob}.

Soit $N\in \{N',N''\}$. On a par \eqref{ineqsupch}
\begin{equation}\label{ineps0}
\log\frac{\si(n)}{n}-\vep\log \tau(n)\le \log\frac{\si(N)}{N}-\vep\log \tau(N)
\end{equation}
et
\begin{equation}\label{ineps}
\log\frac{\si(n)}{n}-\log F(\log\tau(n)) \le g(\log\tau(n))+
\log\frac{\si(N)}{N}-\vep\log \tau(N)
\end{equation}
en posant $g(u)=\vep u-\log F(u)$. La fonction $g$ est convexe sur l'intervalle
$[u_0,u_1] \supset [\log \tau(N'),\log \tau(N'')]$. On choisit $N=N'$ 
ou $N=N''$ pour que
$$g(\log \tau(n))\le \max \left(g(\log \tau(N')),g(\log \tau(N''))\right)=
g(\log \tau(N))$$
et l'inégalité \eqref{ineps} entraîne
\begin{eqnarray*}
\log\frac{\si(n)}{n}-\log F(\log\tau(n)) &\le& g(\log \tau(N))+
\log\frac{\si(N)}{N}-\vep\log \tau(N)\\
&=&\log\;\frac{\si(N)}{N\;F(\log(\tau(N)))}
\end{eqnarray*}
ce qui complète la preuve du lemme \ref{lem4}.
\end{dem}

\begin{lem}\label{lem5}
Soit $a_1,a_2,a_3,a_4,a_5$ des nombres réels vérifiant 
$a1>0,a_2\ge0,a_3\ge 0,a_4\ge 1,a_5\ge e$. La fonction
$$F(u)=a_1\log(a_4+u)+a_2\log\log(a_5+u)+a_3$$
est logarithmiquement concave pour $u > 0$.
\end{lem}

\begin{dem}
On calcule $F'(u)$, $F''(u)$ et l'on montre que, pour $u > 0$, on a
$F(u) > 0$, $F'(u) > 0$ et $F''(u) < 0$ ce qui entraîne 
$F(u)F''(u)-F'^2(u) < 0$.
\end{dem}

\subsection{Les grandes valeurs de $\mathbf{\om(n)}$}

\begin{lem}\label{lem6}
Soit $k_0=15985$ et $n$ tel que $\om(n)=k\ge k_0$. On a
\begin{equation}\label{ineg1}
\frac{\si(n)}{n}\le \frac{n}{\vfi(n)} \le 2.32 \;\log\log\tau(n).
\end{equation}
Soit $n$ tel que $\om(n)=k\ge k_1=166000$, on a
\begin{equation}\label{ineg2}
\frac{\si(n)}{n}\le \frac{n}{\vfi(n)} \le e^\ga \log\log\tau(n)+
e^\ga \log\log\log\tau(n) +0.94\;.
\end{equation}
\end{lem}

\begin{dem}
La démonstration de ce lemme est une forme effective de la preuve de 
\eqref{tenb} dans \cite[formule (17)]{Ten1}.

On définit $\la_k$ par $p_k=k(\log k+\log\log k-\la_k)$ de sorte que, 
par le théorème A de \cite{MasRob}, on a pour $k\ge k_0$,
\begin{equation}\label{lamk}
\la_k =\log k+\log\log k-\frac{p_k}{k} \ge 0.9427.
\end{equation}
Nous utiliserons également la majoration (4.10) de \cite{RS}
et le théorème 6.12 de \cite{Dus}:
\begin{equation}\label{majNk}
\frac{N_k}{\vfi(N_k)}=\prod_{p\;\le\; p_k} \frac{p}{p-1}\le 
e^\ga\left(\log p_k+\de(k)\right) 
\end{equation}
avec
\begin{equation}\label{delta}
\de(k)=
\begin{cases}
\frac{2}{\sqrt{p_k}} &\text{ si } p_k\le 10^8\\
\frac{0.2}{\log p_k} &\text{ si } p_k > 10^8.
\end{cases}
\end{equation}
De la définition de $\la_k$, on déduit pour $k\ge k_0=15985$,
\begin{eqnarray}\label{majpk}
p_k \!\!\!&\le&\!\!\!  k\log k \left( 1+\frac{\log\log k-\la_k}{\log k} 
\right),\notag\\
\log p_k \!\!\!&\le&\!\!\!  \log k+\log\log k +\frac{\log\log k-\la_k}{\log k} 
\le \log k+\log\log k +\be(k)
\end{eqnarray}
avec, par \eqref{lamk}
\begin{equation}\label{beta}
\be(k)=\frac{\log\log k -0.9427}{\log k}\cdot 
\end{equation}
Soit maintenant $n=q_1^{\al_1}q_2^{\al_2}\ldots q_k^{\al_k}$ un 
nombre tel que $\om(n)=k$. On a, comme en \eqref{maj1},
\begin{equation*}
\frac{\si(n)}{n}\le \frac{n}{\vfi(n)}=\prod_{i=1}^k \frac{1}{1-1/q_i}\le 
\prod_{i=1}^k \frac{1}{1-1/p_i}=\frac{N_k}{\vfi(N_k)}  
\end{equation*} 
et, par \eqref{majNk},
\begin{equation}\label{majsi}
\frac{\si(n)}{n}\le \frac{n}{\vfi(n)} \le e^\ga(\log p_k+\de(k)).
\end{equation}
Ainsi, \eqref{majpk} et \eqref{majsi} entraînent pour $k=\om(n)\ge k_0=15985$ 
\begin{equation}\label{majsi2}
\frac{\si(n)}{n}\le \frac{n}{\vfi(n)} \le 
e^\ga(\log k+\log \log k+\be(k)+\de(k)).
\end{equation}

\ni
{\bf {Démonstration de \eqref{ineg1}.}} Comme $p_{k_0}=175939$, on a, 
par \eqref{delta}, pour $k\ge k_0$,
$$\de(k)\le \max\left(\frac{2}{\sqrt{175939}},\;\frac{0.2}{\log(10^8)}\right)=
\frac{0.2}{\log(10^8)}=0.010857\ldots$$
et l'inégalité \eqref{majsi2} implique pour $k\ge k_0$
$$ \frac{n}{\vfi(n)} \le 
e^\ga\log k\left(1+\frac{\log \log k_0+\be(k_0)+0.01086}
{\log k_0}\right) \le 2.23\; \log k.$$
Mais $\tau(n)\ge 2^{\om(n)}=2^k\ge 2^{k_0}$, donc 
\begin{eqnarray*}
\frac{n}{\vfi(n)} &\le& 2.23\; (\log \log \tau(n)-\log\log 2)\\
&\le& 2.23 \;\log \log \tau(n)\left(1+\frac{0.37}{\log\log \tau(n)}\right)\\
&\le& 2.23 \;\log \log \tau(n)\left(1+\frac{0.37}{\log\log 2^{k_0}}\right)\\
&=& 2.31776\ldots \;\log \log \tau(n)
\end{eqnarray*}
ce qui prouve \eqref{ineg1}.

\ni
{\bf {Démonstration de \eqref{ineg2}.}} Puisque $\tau(n)\ge 2^{\om(n)}$, on a
$k=\om(n)\le \frac{\log \tau(n)}{\log 2}$ et il suit
\begin{eqnarray*}
\log k &\le& \log\log\tau(n)-\log\log 2=
\log\log\tau(n)\left(1-\frac{\log\log 2}{\log\log\tau(n)}\right)\\
\log\log k&\le& \log\log\log\tau(n)-\frac{\log\log 2}{\log\log\tau(n)}
\le\log\log\log\tau(n)+\eta(k)
\end{eqnarray*} 
avec $\eta(k)=\frac{-\log\log 2}{\log\log 2^{k}}$ et 
\eqref{majsi2} entraîne
$$\frac{n}{\vfi(n)} \le e^\ga\log \log\tau(n)+e^\ga \log \log \log\tau(n)
+\rho(k)$$
avec
$$\rho(k)=e^\ga\left(-\log\log 2+\eta(k)+\be(k)+\de(k)\right).$$
Pour $k \ge k_2=\pi(10^8)+1=5761456$, on a $p_k \ge p_{k_2}= 10^8+7$,
chacune des fonctions $\eta(k),\be(k)$ et $\de(k)$ est décroissante et 
$$\rho(k) \le \rho(k_2) = 0.921296\ldots$$
Et pour $k1=166000 \le k < k2$, comme $p_{k_1}=2248723$, on a
$$\rho(k) \le \rho(k_1)= 0.939945\ldots$$
ce qui prouve \eqref{ineg2}.
\end{dem}

\subsection{Preuve de \eqref{in3}~: les petites valeurs de $\mathbf{\tau(n)}$}

Soit toujours $k_0=15985,\; p_{k_0}=175939$. On pose 
$\vep^{(0)}=\log\left(1+\frac{1}{p_{k_0}}\right)/\log 2$
 et l'on génère les nombres $(\si,\tau)$--superchampions 
$N$ vérifiant $2\le N\le  N^+_{\vep^{(0)}}$ (cf. définition \ref{defN+}); 
pour chacun d'entre eux, on calcule
\begin{equation*}
f_1(N)=\frac{\si(N)}{N\log\log(3\tau(N))}\cdot
\end{equation*}
Soit $\vep^{(1)}=\log(1+1/113)/\log 2\ap0.012711$; le maximum de $f_1(N)$ 
est atteint pour 
$$N=M_1=N^+_{\vep^{(1)}}=2^83^55^37^211^213^2 \prod_{17\le p \le 113}\; p$$ 
et vaut $2.597907\ldots$ 

On pose $F(u)=\log(u+\log 3)$. Par le lemme \ref{lem5}, $F(u)$ est 
logarithmiquement concave pour $u > 0$. 
 Par la remarque 
\ref{remN+}, il existe $i_0$ tel que $N^+_{(\vep_0)}=N^{(i_0)}$. 
En appliquant le lemme
\ref{lem4} pour tous les couples $(N^{(i)},N^{(i+1)})$ 
vérifiant $2\le N^{(i)} < 
N^{(i+1)}\le N^+_{\vep^{(0)}}=N^{(i_0)}$, on obtient \eqref{in3} 
pour tout $n$ tel que 
$\tau(2)=2\le \tau(n) \le \tau(N^+_{\vep^{(0)}})$.

\subsection{Preuve de \eqref{in3}~: les grandes valeurs 
de $\mathbf{\tau(n)}$}\label{gr1}

Supposons maintenant  $\tau(n) >\tau( N^+_{\vep^{(0)}})=\tau(N^{(i_0)})$.  
Par \eqref{psieps}, $p_{k_0}$ divise $N^+_{(\vep_0)}=N^{(i_0)}$. Soit 
$N^{(i)}$ et $N^{(i+1)}$ les deux nombres 
$(\si,\tau)$--superchampions tels que
$\tau(N^{(i)})\le \tau(n) < \tau(N^{(i+1)})$
(cf. lemme \ref{lemsupch} et remarque \ref{remN+}). 
On doit avoir $i\ge i_0$ et, par le lemme \ref{lemexpo}, 
les nombres premiers $p\le p_{k_0}$ divisent $N^{(i)}$ et $N^{(i+1)}$; ainsi
$\om(N^{(i)})\ge k_0$ et $\om(N^{(i+1)})\ge k_0$. Par \eqref{ineg1}, on a alors
$f_1(N^{(i)})\le 2.32$,  $f_1(N^{(i+1)})\le 2.32$, ce qui entraîne 
par le lemme 
\ref{lem4}, $f_1(n)\le 2.32$ et achève la preuve de \eqref{in3}.

\subsection{Preuve de \eqref{in4}~: les petites valeurs de $\mathbf{\tau(n)}$}

La preuve de \eqref{in4} est très voisine. On pose
$k_1=166000,\; p_{k_1}=2248723$, 
$\vep^{(2)}=\log(1+{1}/{p_{k_1}})/\log 2 \ap 0.00000064156$
 et pour chaque nombre $(\si,\tau)$--superchampion 
$N$ vérifiant $2\le N\le  N^+_{\vep^{(2)}}$, on calcule 
$$f_2(N)=\frac{\si(N)}{N}-e^\ga\log\log(e\tau(N))-
e^\ga\log\log\log(e^e\tau(N)).$$

Soit $\vep^{(3)}=\log(1+1/45439)/\log 2\approx 0.0000317498$; le maximum 
de $f_2(N)$ pour $N\le N^+_{\vep^{(2)}}$ est atteint pour 
$$N=M_2=N^+_{\vep^{(3)}}=2^{18}3^{11}5^77^6
\prod_{11\le p \le 19}\; p^4\prod_{23\le p \le 47}\; p^3
\prod_{53\le p \le 277}\; p^2 \prod_{281\le p \le 45439}\; p$$ 
et vaut
$\mu=f_2(M_2)\approx 0.9414440795$.

Par le lemme \ref{lem5}, la fonction 
\begin{equation}\label{F}
F(u)=e^\ga\log(1+u)+e^\ga\log\log(e+u)+\mu
\end{equation}
est logarithmiquement concave.
D'après notre calcul, pour chaque nombre $(\si,\!\tau)$--superchampion,
$2\le N\le N^+_{\vep^{(2)}}$, on a $\frac{\si(N)}{NF(\log\tau(N))}\le 1$ et 
par le lemme \ref{lem4}, cela entraîne 
que pour tout $n$ vérifiant $\tau(2)=2\le \tau(n)\le \tau(N^+_{\vep^{(2)}})$, 
on a $\frac{\si(n)}{nF(\log\tau(n))}\le 1$,
c'est-à-dire \eqref{in4}.

\subsection{Preuve de \eqref{in4}~: les grandes valeurs de $\mathbf{\tau(n)}$}

Soit maintenant $N$ un nombre $(\si,\tau)$--superchampion supérieur ou égal à 
$N^+_{\vep^{(2)}}$. Comme dans le paragraphe \ref{gr1}, on a $\om(N)\ge k_1$ 
et donc, par \eqref{ineg2} et \eqref{F},
\begin{equation}\label{FN}
\frac{\si(N)}{N} \le e^\ga \log\log \tau(N)+e^\ga \log\log\log \tau(N)+ 0.94
< F(\log\tau(N)).
\end{equation}
Supposons $\tau(n) > \tau(N^+_{\vep^{(2)}})$. On définit comme en \ref{gr1}
les deux nombres $(\si,\tau)$--superchampions $N^{(i)}$ et $N^{(i+1)}$ tels
que $\tau(N^{(i)}) \le \tau(n) < \tau(N^{(i+1)})$. On a $\om(N^{(i+1)})
\ge \om(N^{(i)})\ge \om(N^+_{\vep^{(2)}})\ge k_1$. Par \eqref{FN}, il vient
$$\max_{N\in \{N^{(i)},N^{(i+1)}\}} \frac{\si(N)}{NF(\log\tau(N))} < 1$$
d'où l'on déduit, par le lemme  \ref{lem4}, que
$\frac{\si(n)}{n} < F(\log\tau(n))$, et \eqref{in4} est démontrée.

\mk

Notons que, pour $N=M_2$, les valeurs de $\frac{\si(N)}{N}\approx
19.0983$ et $\frac{N}{\vfi(N)}\approx 19.1096$ sont très voisines.

\section{Démonstration du théorème \ref{th3}}\label{par5}

\subsection{La méthode des bénéfices}

\begin{definition}
Soit $\vep > 0$ et soit $N$ l'un des nombres $(\si,\tau)$--superchampions
qui maximisent la fonction $n\mapsto \frac{\si(n)}{n\tau(n)^\vep}$. Pour
$n\ge 1$, on appelle {\it bénéfice} de $n$ la quantité (qui dépend de $\vep$)
\begin{equation}\label{defben}
\ben (n)=\log\frac{\si(N)}{N\tau(N)^\vep}-\log\frac{\si(n)}{n\tau(n)^\vep}=
\log \frac{\si(N)/N}{\si(n)/n}-\vep\log \frac{\tau(N)}{\tau(n)}\cdot
\end{equation}
\end{definition}
Compte tenu de \eqref{ineqsupch}, on a
\begin{equation}\label{benpos}
\ben (n) \ge 0.
\end{equation}
Soit $n=\prod_{p} \; p^{\be_p}$ et $N=\prod_{p} \; p^{\al_p}$ (où pour 
chaque $p$ premier, $\al_p$ vérifie \eqref{psieps}). On a
\begin{equation}\label{benbenp}
\ben (n)=\sum_{p\;\text{ premier }} \; \benp (n)
\end{equation}
avec
\begin{eqnarray}\label{benp}
\benp (n) &=&
\log\;\frac{\si(p^{\al_p})/p^{\al_p}}{\si(p^{\be_p})/p^{\be_p}}
-\vep\log\;\frac{\tau(p^{\al_p})}{\tau(p^{\be_p})}\notag\\
&=& \log\left(\frac{1+\frac{1}{p}+\ldots+\frac{1}{p^{\al_p}}}
{1+\frac{1}{p}+\ldots+\frac{1}{p^{\be_p}}}\right)-
\vep\log\left(\frac{\al_p+1}{\be_p+1}\right)\notag\\
&=& \log\left(\frac{1-\frac{1}{p^{\al_p+1}}}
{1-\frac{1}{p^{\be_p+1}}}\right)-
\vep\log\left(\frac{\al_p+1}{\be_p+1}\right)\cdot
\end{eqnarray}
\begin{lem}\label{lemben}
(i) Les nombres $\al_p$ et $\vep$ étant liés par \eqref{psieps}, 
la quantité $\benp (n)$ définie par \eqref{benp} vérifie $\benp(n)\ge 0$ pour
tout $\be_p\ge 0$.
De plus, la fonction $\be_p\mapsto \benp (n)$ est 
croissante (au sens large) pour $\be_p\ge \al_p$, 
décroissante (au sens large) pour $\be_p\le \al_p$, 
nulle pour $\be_p=\al_p$ et tend vers l'infini avec $\be_p$.

\ni
(ii) Pour $\vep$ fixé et $\psi(p,1) < \vep$, par \eqref{psieps} on a $\al_p=0$.
Lorsque $\be_p=1$, $\benp (n)=\vep\log 2-\log(1+1/p)$ est une fonction 
croissante en $p$ qui tend vers $\vep\log 2$ quand $p\to\iy$. 
\end{lem}

\begin{dem}
Supposons $\be_p\ge\al_p$. A l'aide des fonctions $\psi$ et $\te$ définies en 
\eqref{psi} et en \eqref{te}, on a par \eqref{tem}
\begin{eqnarray}\label{benp1}
\benp (n)=\te(p,\al_p)-\te(p,\be_p) &=&
\sum_{a=\al_p+1}^{\be_p} \te(p,a-1)-\te(p,a)\notag\\
&=& \sum_{a=\al_p+1}^{\be_p} \log\left(1+\frac 1a\right) 
\left[ \vep-\psi(p,a)\right]
\end{eqnarray}
ce qui, par la décroissance de la fonction $\psi$ (cf. lemme \eqref{lempsi}) et
\eqref{psieps} prouve $\benp(n)\ge 0$ et 
la croissance de $\benp(n)$ par rapport à $\be_p$.

Si $0\le \be_p < \al_p$, on a
\begin{equation}\label{benp2}
\benp (n)=\sum_{a=\be_p+1}^{\al_p} \log\left(1+\frac 1a\right) 
\left[\psi(p,a)-\vep\right]
\end{equation}
et l'on conclut de la même façon pour montrer
la décroissance de $\benp(n)$ par rapport à $\be_p$.

Lorsque $p,\vep$ et $\al_p$ sont fixés, la formule \eqref{benp} montre que 
$\benp(n)$ tend vers l'infini avec $\be_p$.

La preuve de (ii) est facile.
\end{dem}

\begin{coro}\label{corben}
Soit $\vep$ fixé et $B$ un nombre réel, $0\le B < \vep\log 2$.
L'ensemble des nombres entiers $n$ vérifiant $\ben(n)\le B$ (où $\ben(n)$
est défini par \eqref{defben}) est fini.
\end{coro}

\begin{dem}
Soit $n$ un nombre vérifiant $\ben(n)\le B$.
Par le lemme \ref{lemben} (i), on a $\benp(n)\ge 0$ et par \eqref{benbenp},
on a, pour chaque nombre premier $p$, $\benp(n)\le \ben(n)\le B < \vep\log 2$.
Toujours par le lemme \ref{lemben}, cela implique, lorsque $\al_p=0$, 
qu'il existe $p_0$ tel que $\be_p=0$ pour $p > p_0$. 

Ainsi, les nombres premiers supérieurs à $p_0$ ne divisent pas $n$. 

Soit maintenant $p\le p_0$. Le lemme \ref{lemben} (i) et l'hypothèse
$\benp(n)\le B$ montrent qu'il n'y a qu'un nombre fini d'exposants $\be_p$
possibles.  
\end{dem}

\begin{rem}\label{remben}
La démonstration du corollaire \ref{corben} est effective et permet,
lorsque $B$ est petit, de déterminer l'ensemble $\cn(B)$ des nombres $n$
vérifiant $\ben(n)\le B$.

Cependant il semble très difficile d'obtenir une estimation de Card$(\cn(B))$
en fonction de $B$.
\end{rem}

Nous pouvons maintenant préciser le lemme \ref{lem4}.

\begin{lem}\label{lem4bis}
Si l'on rajoute dans le lemme \ref{lem4} la condition
\begin{equation}\label{ben<=B}
\ben(n)\ge B 
\end{equation}
alors la conclusion devient
$$\frac{\si(n)}{n\;F(\log(\tau(n)))}\le e^{-B} \;\max
\left(\frac{\si(N')}{N'\;F(\log(\tau(N')))}\;,\;\;
\frac{\si(N'')}{N''\;F(\log(\tau(N'')))}\right).$$
\end{lem}

\begin{dem}
Compte tenu de \eqref{defben} et de \eqref{ben<=B}, l'inégalité
\eqref{ineps0} se réécrit :
\begin{eqnarray}\label{ineps2}
\log\frac{\si(n)}{n}-\vep\log \tau(n) &=& 
\log\frac{\si(N)}{N}-\vep\log \tau(N)-\ben(n)\notag\\
&\le& \log\frac{\si(N)}{N}-\vep\log \tau(N)-B.
\end{eqnarray}
Ainsi, \eqref{ineps} devient
\begin{equation*}\label{ineps3}
\log\frac{\si(n)}{n}-\log F(\log\tau(n)) \le g(\log \tau(n))+
\log\frac{\si(N)}{N}-\vep\log \tau(N)-B
\end{equation*}
et la démonstration se termine comme celle du lemme \ref{lem4}.
\end{dem}

\medskip
\ni
{\bf Démonstration du théorème \ref{th3}~}: Les calculs effectués pour 
 démontrer la formule \eqref{in3} du théorème \ref{th2} 
montrent que les seuls nombres $(\si,\tau)$--superchampions $N$ tels que
$$f_1(N)=\frac{\si(N)}{N\log\log(3\tau(N))} \ge \frac{2957}{1000}$$
sont les nombres $N=N^{(i)}$ (cf. lemme \ref{lemsupch}) avec $46\le i\le 50$.
La table ci-dessous donne la valeur de ces nombres en fonction de
$$M_1=N^{(46)}=2^83^55^37^211^213^2 \prod_{17\le p \le 113}\; p\;.$$ 
Notons que pour $45\le i\le 51$, $\vep_i=\psi(p^{(i)},1)=
\log(1+1/p^{(i)})/\log 2$.

$$
\begin{array}{|r|c|c|r|c|}
i&\vep_i&p^{(i)}&N^{(i)}/M_1&f_1(N^{(i)})\\
\hline
45&0.0132&109&1/113&2.596216\\
46&0.0127&113&1&2.597907\\
47&0.0113&127&127&2.597801\\
48&0.0110&131&127\times131&2.597746\\
49&0.0105&137&127\times131\times137&2.597461\\
50&0.0103&139&127\times131\times137\times139&2.597502\\
51&0.0097&149&127\times131\times137\times139\times149&2.596862\\
\hline
\end{array}
$$

\mk
\ni
Ainsi, par le lemme \ref{lem4},
\eqref{in3bis} est satisfaite pour tous les nombres $n$ vérifiant
$2\le \tau(n)\le \tau(N^{(45)})=\frac 12\;\tau(M_1)$ ou $\tau(n)\ge 
\tau(N^{51)})=32\;\tau(M_1)$.

Ensuite, on pose $B=\log \left(\frac{2.6}{2.597}\right)\approx 0.0011545$ 
de telle sorte que, 
par \eqref{in3}, pour tout $n\ge 2$, on ait
$$ e^{-B}f_1(n)=e^{-B}\frac{\si(n)}{n\log\log(3\tau(n))} 
< e^{-B}\;\times \; 2.6=\frac{2957}{1000}$$
ce qui, par le lemme \ref{lem4bis},  prouvera \eqref{in3bis} pour les $n$ 
restants qui ont un bénéfice supérieur à $B$.

Finalement, pour chacune des 6 valeurs de $\vep=\vep_i$ avec $45\le i\le 50$,
on détermine les $n$ dont le bénéfice est inférieur à $B$
(cf. corollaire \ref{corben}). Parmi ces nombres
$n$, seuls 12 présentent une valeur de $f_1(n)$ supérieure à $2.597$. Ils sont 
énumérés par valeur décroissante de $f_1(n)$ dans le tableau ci-dessous.

$$
\begin{array}{|r|c|c|}
n/M_1\;&\tau(n)/\tau(M_1)&\quad f_1(n)\quad\\
\hline
1\; &1&2.597907\\
127\; &2&2.597801\\
127\times131\; &4&2.597746\\
127\times131\times137\times139\; &16&2.597502\\
127\times131\times137\; &8&2.597461\\
 2\times127\times131\; &40/9&2.597331\\
2\; &10/9&2.597290\\
2\times 127\; &20/9&2.597288\\
\quad2\times127\times131\times137\times139\; &160/9&2.597269\\
 127\times131\times139\; &8&2.597190\\
131\; &2&2.597181\\
2\times127\times131\times137\; &80/9&2.597140\\
\hline
\end{array}
$$

\bk

Si l'on désigne par $\nu(x)$ le cardinal de l'ensemble des nombres $n$ pour 
lesquels $f_1(n)=\frac{\si(n)}{n\log\log \tau(n)}\ge x$, nous avons 
calculé les valeurs suivantes.

$$
\begin{array}{|r|cccccccc|}
x&2.597&2.596&2.595&2.594&2.593&2.592&2.591&2.590\\
\hline
\nu(x)&12&45&179&586&1680&4760&12653&32187\\
\hline
\end{array}
$$

\bk

Posons
$$ M_3 = 2^{9}3^{5}5^37^311^213^217^2\prod_{19\le p \le 211}\; p$$
et considérons la suite des nombres $n_p=\frac{pM_3}{211}$ où $p$ parcourt 
les nombres premiers supérieurs à $211$. On a $\tau(n_p)=\tau(M_3)$,
$f_1(n_p)=\frac{211(p+1)}{212\; p}f_1(M_3)$ et 
$$\lim_{p\to\iy} f_1(n_p)
= \frac{211}{212}f_1(M_3)=2.580303\ldots$$
Ainsi, pour $x < 2.58$, la quantité $\nu(x)$ est infinie.

Il est vraisemblable que l'on a
$$\limsup_{n\to\iy} \frac{\si(n)}{n\log \log (3\tau(n))}=
\frac{211}{212}\frac{\si(M_3)}{M_3\log \log (3\tau(M_3))}=
2.580303\ldots$$

\subsection*{Remerciements.} J'ai plaisir à remercier A. Schinzel 
pour la question qu'il m'a posée, et qui m'a conduit au théorème \ref{th1}.
Je remercie également M. Deléglise pour le dessin de la figure et pour 
d'utiles discussions au sujet des nombres superchampions.

\def\refname{Références}

\vspace{2cm}

\noindent
Jean-Louis Nicolas,

\noindent
Universit\'e de Lyon, Universit\'e de Lyon 1, CNRS,

\noindent
Institut Camille Jordan, UMR 5208, 

\noindent
B\^at. Doyen Jean Braconnier,

\noindent
21 Avenue Claude Bernard,

\noindent
F-69622 Villeurbanne c\'edex, France.

\bigskip
    
\ni
\texttt{jlnicola@in2p3.fr} 

\ni
\texttt{http://math.univ-lyon1.fr/$\sim$nicolas/}.

\end{document}